# A Dynamic Theory-Based Method for the Computation of an Unstable Equilibrium Point

Robert Owusu-Mireku, *Student Member,* Matt Hin, *Student Member, IEEE*, Hsiao-Dong Chiang, *Fellow, IEEE*

*Abstract*—In this paper, a new combination of a dynamic transformation method and a trajectory-based integration technique is proposed for the model independent computation of unstable equilibrium points (UEPs). The transformation method converts a UEP into a stable equilibrium point (SEP) to expand the convergence region by creating a quotient gradient system. The resulting SEP is then calculated using a quasi–Newton form of the pseudo-transient continuation method that exploits the structure of the quotient gradient system to speed up computation. The proposed method's conditions for convergence are presented, and the method is tested on the WSCC 9-bus 3-machine system and the IEEE 145-bus 50-machine system. The results show that the proposed method gives accurate results, it is sufficiently fast, numerically stable, and enlarges the convergence region of the UEP.

## I. INTRODUCTION

THE fast and accurate computation of unstable equilibrium points (UEPs) is important in most direct method applications in power systems [1, 2, 5-8, 34]. In the closest and controlling UEP methods for transient stability assessment, the computation of type-1 UEPs is required to determine the relevant critical energy level needed for stability assessment [1, 2, 8, 34]. The computation of a UEP is also required in assessing the proximity of an operating point to a system's voltage collapse limit [5-7] for voltage stability assessment [3]. UEPs are also computed in power flow analysis [4] when all the possible solutions of a power system are required [9–13]. However, the computation of UEPs, with algebraic solvers like Newton-Raphson's method, are inherently difficult due to the generally small size of their convergence regions [1, 30] and the difficulty in finding an initial point that is sufficiently close to the UEP.

Various methods have been proposed for computing UEPs. Some methods like the corrected corner, the ray point approximation method, the MOD, the BCU, and the shadowing methods [1, 14-17, 26] are based on finding efficient and robust ways of determining the correct initial point for the UEP, and then solving for the exact UEP using algebraic solvers like the Newton-Raphson method. Other methods are based on continuation and homotopy methods [9–13, 18, 20, 22, 23] or optimization techniques [27].

Even though some of these methods, like the BCU method, have theoretical bases they have been found to fail sometimes [1, 14-17]. Others like the homotopy and optimization techniques also tend to be intractable in large system applications [9–13, 18, 20, 22, 23, 27].

In [19, 21, 24], the authors addressed the challenge that accompanies computing UEPs with the above methods by proposing methods that transformed the problem into the computation of stable equilibrium points (SEPs), an approach generally called dynamic methods. The methods proposed in [21, 24] used a spectral decomposition of the Jacobian of the original system to construct a new gradient system where the UEPs in the original system are SEPs, while [19] used the whole Jacobian of the original system for the construction of the new gradient system. The proposed methods in [21, 24] are only applicable to the power system network-reduction models, and the method in [19] focuses more on using the transformation and related minus systems to determine initial guesses sufficiently close to the UEPs on the stability boundary of an SEP. The author in [25] has shown that generally, transformation techniques like the one in [21] do not work properly when used to compute the closest UEPs on a stability boundary.

In this paper, we propose a new theory-based solver for the robust computation of a UEP of a system of ordinary differential equations (ODEs) or differential algebraic equations (DAE) given an initial guess. It combines two steps: 1) the transformation of the (relevant UEP of a given system into a stable equilibrium point (SEP) of a related dynamical system, called the quotient gradient system (QGS) [28]; and (2) computation of the resulting SEP of the QGS using a trajectory-unified method (TJU) like the exact or the in-exact/quasi Pseudo-transient continuation ($\psi$tc) method [33]. Numerical results suggest that the proposed method generally has better convergence regions as compared to the Newton-Raphson (NR) method when applied to UEP computation but sacrifices some speed for the robustness. The proposed method is to be used in tandem with existing methods, like the BCU method, or any of the other methods that require an algebraic solver for UEP computations.

The QGS transformation converts the UEP into an SEP, which has a stability region, and hence, can be solved with a

Robert Owusu-Mireku and Hsiao-Dong Chiang are with the Department of Electrical and Computer Engineering, Cornell University, Ithaca, NY 14853, USA (e-mail: ro82@cornell.edu, hc63@cornell.edu)

Matt Hin is with the Center Applied Mathematics, Cornell University, Ithaca, NY 14853, USA (e-mail: mfh72@cornell.edu).

TJU method. Also, pairing the transformation with $\psi$tc expands the resulting convergence region of the UEP, since the $\psi$tc has a more continuous (connected and smooth) convergence region, unlike NR, where the convergence region is a fractal [29].

The major difference between our proposed method and the method in [19, 21, 24] is that the proposed method does not require eigenvalue computations, it uses a fast QGS structure dependent quasi-Newton method to solve for the UEP after the transformation, and it is applicable to UEP computations in general. The transformation method used in [19] is similar to the transformation method used in our work, but the authors in [19] focus on determining a initial points close to UEPs, and also exploiting the periodic structure and Lyapunov function of the resulting system to validate the calculated UEPs. [19] also focuses on finding multiple UEPs while this work focuses on finding a specific UEP given an initial condition. In this work, we also focus on exploiting the structure of the QGS resulting from the transformation to speed up the computation of the UEP and to ensure convergence to a solution of the original system. The use of this specific QGS structure dependent in-exact/quasi method to speed up computation and avoid convergence to points that are not solutions to the original system is new in this kind of power system application. In effect, this proposed method will provide a theory-based method to improve the robustness of UEP computations, by converting the UEP into an SEP and enlarging its convergence region. The method is also independent of the network model and can be applied to problems that are unrelated to power systems. Our proposed method could be used in tandem with [19] for solving for multiple UEPs in a power system. Under certain assumptions, the proposed method can have local q-superlinear or local quadratic convergence.

This paper is organized as follows. Section II presents the transformation of the original system. Section III reviews the pseudo-transient continuation ($\psi tc$) method for ODE systems. Section IV proposes the QGS-based pseudo-transient continuation ($\psi tc$) method and presents the conditions for its convergence. Section V presents two numerical examples, discusses the findings from the simulations, and proposes an algorithm that combines the proposed method with the NR method for computing UEPs. The conclusions are then stated in Section VI.

## II. SYSTEM TRANSFORMATION

### A. Original Problem Formulation

It is a known fact that the equilibrium points of an ODE or DAE can be found by solving for the zeros of its corresponding equilibrium equations. Without loss of generality, the equilibrium equation of an ODE or DAE can be represented by (1).

$$F(x) = 0 \qquad (1)$$

where $F: R^n \rightarrow R^n$, $F$ is assumed to be $C^2$, $x \in R^n$ is a vector of equilibrium states, and $n \geq 1$.

$F(x)$ is either the vector field of the system of ODEs or the vector field and the algebraic manifold of the system DAEs. Regardless of what (1) represents, UEPs are inherently very difficult to compute because appropriate initial guesses are difficult to determine, and UEPs generally have a small convergence region with respect to a numerical method, say, Newton-Raphson [1, 29]. However, like any equilibrium point or zero of a function, UEPs can be computed using algebraic solvers when an initial point sufficiently close to the UEP is provided. The size and continuity/compactness of the convergence region of an algebraic solver determines how close an initial guess must be for the solver to successfully converge to the UEP. Methods like the BCU boundary following procedure [1] has approached this problem by focusing on finding an initial point sufficiently close to the UEP, but sometimes finding an initial point close enough to the UEP for an algebraic solver like NR to work is not possible. The purpose of this work is to provide a fast solver that has a large connected convergence region, implying that the initial guesses can be further away from a UEP. One common solver used for UEP computations is the NR method.

$$\dot{x} = f(x) = -DF(x)^{-1}F(x) \qquad (2)$$

Each step in the NR algorithm can also be viewed as a forward Euler step of the dynamic system (2) with a time step of 1, where $DF(x)$ is the Jacobian matrix of $F(x)$ and $-DF(x)^{-1}F(x)$ is the vector field $f(x)$ of the ODE (2) [30]. Thus, the NR method and some of its variants—for example, Iwamoto's method—are basically the forward Euler integration of the new dynamic system (2) [30]. This new system (2) is stable at all equilibrium points where $DF(x)$ is nonsingular. Consequently, in terms of the computation of UEPs, the NR method and some of its variants can be considered as a numerical technique that involves the transformation of a UEP of the original system (1) into an SEP of (2), and an application of an explicit integration method like the Euler method, to solve for the new system's SEP, which is the UEP of the system of differential equations corresponding to (1) for a given initial point.

Unlike the transformation used in this work, the transformation used in the NR method employs the inverse of a Jacobian $DF(x)$ which can be singular somewhere in the neighborhood of the resulting SEP. Implying that (2) does not always satisfy the requirements for existence and uniqueness of solution.

Some variants of the NR method, like the Continuous Newton-Raphson method (Continuous NR) presented in [30], propose the use of a much more stable explicit integration technique, such as the fourth order Runge-Kutta (RK4) method over the forward Euler method, for the integration of the dynamic system in (2).

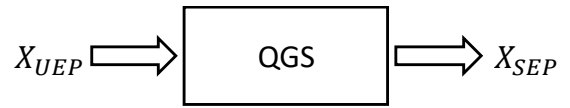

Fig. 1. An illustration of the Quotient Gradient System Transformation.

## B. Quotient Gradient System Transformation

The transformation step used in the proposed method is illustrated in Fig. 1. It is called the Quotient Gradient System (QGS) [28] transformation, and has the following form:

$$\dot{x} = Q(x) = -DF(x)^T F(x). \quad (3)$$

*Proposition 1:* Let all the equilibrium points of (3) be hyperbolic and finite in number. If $\hat{x}$ is the solution of (1), then $\hat{x}$ is an SEP of (3).

Proof:

It is obvious that if $\hat{x}$ is a solution of (1), then it is an equilibrium point of (3). For $\hat{x}$ to be an SEP of (3), the Jacobian ($DQ(x)$) of (3) at the equilibrium point $\hat{x}$ should be negative definite.

Thus, for any nonzero vector, $y \in R^n$, $y^T DQ(\hat{x})y < 0$. The Jacobian of (3) at a point $x$ is given by (4):

$$DQ(x) = -\sum_{k=1}^{n} f_k H^T(f_k) - DF(x)^T DF(x) \quad (4)$$

where $f_k$ is the kth function in $F(\hat{x})$ and $H^T(f_k)$ is the Hessian of $f_k$. At $\hat{x}$, $f_k = 0$ and hence,

$$DQ(\hat{x}) = -DF(\hat{x})^T DF(\hat{x}). \quad (5)$$

The Jacobian matrix (5) is negative semi-definite. Since all the equilibrium points of (3) are assumed to be hyperbolic, it implies that none of the equilibrium points have eigenvalues with a zero real part, and hence, none of the equilibrium points have zero eigenvalues. That is, all the eigenvalues of the equilibrium points of (3) have negative real parts, which implies that (5) is negative definite.

The advantage of the proposed method over the NR method is that the stability region of the SEP resulting from the transformation exist and is compact/continuous (connected) if the underlying assumptions are satisfied, while there is no stability region for the SEP from the NR transformation (since the resulting system does not have unique solutions). The UEP and consequently, its corresponding stability region in the new system (3) can also be efficiently computed by solving for the SEP using implicit integration techniques, which are more stable compared to explicit ones like the Euler method.

## III. THE PSEUDO-TRANSIENT CONTINUATION METHOD

### A. Introduction

After transforming the algebraic problem into computation of an SEP of a dynamic system (3), the most efficient way to solve for the equilibrium point is to use the TJU method. Since we are only interested in the steady state solution, we need a TJU method that will converge quickly to the right steady state solution. The TJU method can be explicit (Euler or Runge-Kutta) or implicit (the trapezoidal method). One such implicit method is the pseudo-transient continuation ($\psi tc$) method. The $\psi tc$ is an implicit TJU that employs adaptive time-stepping for the computation of steady state solutions for partial differential equations, ODEs, and semi-explicit index-one DAEs [32-33]. It is analogous to an implicit integration of a dynamic system with increasing time steps as the system trajectory approaches the steady state solution.

The trajectory-based nature of $\psi tc$ makes its convergence region a better approximation of an SEP's stability region compared to the convergence region of other algebraic solvers. This characteristic of $\psi tc$ implies that the initial points don't have to be as close as the NR method requires. It also means trajectories are more likely to converge to the correct physical SEP and not to other SEPs or non-physical local minima [33]. The adaptive time steps used in the $\psi tc$ method makes it faster than conventional fixed-step integration methods, since larger time steps are taken as the trajectory gets closer to the EP. The $\psi tc$ method is an implicit method, making it numerically more stable than explicit methods like the Euler and Runge-Kutta methods.

Given an initial value problem of the form (6):

$$\dot{x} = -G(x), \quad x(0) = x_0 \quad (6)$$

The steady state solution can be found by integrating (6) with $\psi tc$. Each step in the $\psi tc$ method is given by (7):

$$x_{i+1} = x_i - \left(h_i^{-1} I + DG(x_i)\right)^{-1} G(x_i) \quad (7)$$

where $I$ is an identity matrix of appropriate size, and $h_i$ is a variable time step systematically adjusted to improve rate of convergence to steady state. The time step can be adjusted using the "switch evolution relaxation" (SER) (8) or the norm of the steps (9) [33]:

$$h_i = min\left(h_{i-1} \frac{\|G(x_{i-1})\|}{\|G(x_i)\|}, h_{max}\right) \quad (8)$$

$$h_i = min(h_{i-1} \|x_i - x_{i-1}\|, h_{max}) \quad (9)$$

where $h_{max}$ is a large upper bound of $h_i$. For the results in this work, equation (8) was used for the time step adjustments and $h_{max}$ was set to $\infty$. $h$ can be a vector of different time steps if the system of equations in (6) is stiff or has an ill-conditioned Jacobian.

$\psi tc$ Algorithm [33]:
1. Set $x = x_0$ and $h = h_0$. Evaluate $G(x)$.
2. While $\|G(x)\|$ is larger than a threshold:
   a. Solve $(h^{-1} I + DG(x))s = -G(x)$.
   b. Set $x = x + s$.
   c. Evaluate $G(x)$.
   d. Update $h$.

In [33], the authors prove that if a steady state solution exists, then $\psi tc$ for ODEs of the form (6) has a local q-superlinear or quadratic convergence if some assumptions are satisfied. To improve the computational performance of the application of the $\psi tc$ to (3) and guarantee convergence to solutions of (1), we proposed inexact approach, which we call the Quotient Gradient System-based Pseudo-Transient Continuation method (QGS-based $\psi tc$) and is presented in the next section.

## IV. THE QUOTIENT GRADIENT SYSTEM-BASED PSEUDO-TRANSIENT CONTINUATION METHOD

### A. Proposed Method

If we apply the $\psi tc$ method to (3), each step of the $\psi tc$ method, (7), can be rewritten, as shown in (10). We can then simply solve (10) at step (2a) for each iteration of the $\psi tc$ algorithm without any further modifications, since the UEP is now an SEP due to the QGS transformation. However, this approach requires the construction of two Jacobians at each iteration, one for the QGS transformation and another for step (2a) of the $\psi tc$ algorithm. This can be computationally expensive if the Jacobians are constructed numerically. Also, the analytical Jacobian for step (2a) can be complex and error-prone. The use of automatic differentiation for a Jacobian construction in step (2a) can also be quite challenging for complex functions or systems of equations like (3).

$$x_{i+1} = x_i - \left( h_i^{-1} I + \sum_{k=1}^{n} f_k H^T(f_k) + DF(x_i)^T DF(x_i) \right)^{-1} DF(x_i)^T F(x_i) \quad (10)$$

If we assume that $\|-\sum_{k=1}^{n} f_k H^T(f_k)\|$ is sufficiently small, which is true as we approach the equilibrium point of (3), then we can use a quasi-Newton method approach and approximate the Jacobian of (3) with (11):

$$D\tilde{Q}(x) \approx -DF(x)^T DF(x) \quad (11)$$

$$\left( h^{-1} I + DF(x)^T DF(x) \right) s \approx -DF(x)^T F(x) \quad (12)$$

$$\min_s \left\| \begin{bmatrix} DF(x) \\ h^{-\frac{1}{2}} I \end{bmatrix} s + \begin{bmatrix} F(x) \\ 0 \end{bmatrix} \right\| \quad (13)$$

$$\|(h^{-1}I + DF(x)^T DF(x))s + DF(x)^T F(x)\| \leq \xi \|DF(x)^T F(x)\| \quad (14)$$

$$x_{i+1} = x_i + s. \quad (15)$$

Step (2a) in the $\psi tc$ algorithm can then be replaced by either (12) or (13). For small dense systems, it might be more efficient to solve (12) using QR factorization. In larger sparse systems, solving (13) using Cholesky decomposition or iterative methods like the conjugate gradient methods (precondition conjugate gradient methods) might be more efficient. Generically, (12) and (13) can be represented by (14) where $\xi$ could be related to the difference between the exact Jacobian of (3) and the approximate Jacobian (11) if (12) or (13) are solved using QR factorization or Cholesky decomposition, or $\xi$ could be related to the Jacobian approximation and the inexact steps involved when iterative methods are used to solve the Newton step.

The proposed QGS-based $\psi tc$ method can be summarized by the following steps.
1. Transform the original algebraic problem (1) into a dynamic system (3) using the QGS transformation.
2. Starting at the given initial guess, apply the $\psi tc$ method to the original system's (1) surrogate QGS system, (3), solving either (12) or (13) at step (2a) of the $\psi tc$ algorithm.

In comparison to the NR method, the QGS-based $\psi tc$ method will, in most cases, require more iterations since it does not always converge quadratically to the equilibrium point. Since it will also require more computations per iteration compared to the NR method, we propose that the QGS-based $\psi tc$ be used as a re-starting algorithm after the NR method fails. The idea is to re-start the UEP computation with the QGS-based $\psi tc$ method and then switch back to the NR method when the convergence criterion for the QGS-based $\psi tc$ method is below a defined threshold. It should be noted that the QGS-based $\psi tc$ method shares some similarities with the trust region method used for nonlinear least square problems, but the QGS-based $\psi tc$ method was derived independently.

### B. Convergence

In this section we present the type of convergence and the conditions necessary for the convergence of the QGS-based $\psi tc$ method. We show that the QGS-based $\psi tc$ can have local q-superlinear or even local quadratic convergence if certain assumptions and conditions are satisfied. We also show that the Jacobian approximation used in the QGS-based $\psi tc$ method guarantees that if the QGS-based $\psi tc$ method convergences it will only converge to solutions of the original system (1). Let $x^*$ be a UEP of the dynamic system with equilibrium equations represented by equation (1).

Assumptions:
1. $DF(x)^T F(x)$ is everywhere defined and Lipschitz continuously Fréchet differentiable, $\|DF(x)^T F(x)\| \leq M, M > 0$ for all $x$.
2. There are $\epsilon_2, \beta > 0$ such that if $\|x - x^*\| < \epsilon_2$, then $\left\| \left( h^{-1}I + D(DF(x)^T F(x)) \right)^{-1} \right\| \leq (1 + \beta h)^{-1}$ for all $h \geq 0$.
3. Equation (3) has unique hyperbolic equilibrium points.
4. $\underset{i}{\text{Inf}}\, h_i > 0$. This assumption must be satisfied to prevent the QGS-based $\psi tc$ method from stalling.

*Proposition 2:*
Let $F(x^*) = 0$. Let Assumptions 1,2,3 and 4 hold. The QGS-based $\psi tc$ method will have a local q-superlinear convergence to $x^*$ from an initial point $x_0$ if $x^*$ is a stable hyperbolic equilibrium point of the QGS system (3), $x_0$ is in the stability region of $x^*$ for the QGS system (3), $\{h_i\}$ is of the form (8), $h_{max} = \infty$, $\xi_i \leq \hat{\xi}$, $h_i \to h_{max}$, and $\xi_i \to 0$. The convergence to $x^*$ will be q-linear if $h_{max} < \infty$.

Proof Sketch:
From proposition 1, if $F(x^*) = 0$ and Assumption 3 holds, then $x^*$ is a hyperbolic stable equilibrium point of (3). Consequently, $x^*$ has a stability region in (3). We also know that the portion of $\xi_i$ that is due to the approximation of the Jacobian in the QGS-based $\psi tc$ method is directly proportional to $\|-\sum_{k=1}^{n} f_k H^T(f_k)\|_i$, and that $\|-\sum_{k=1}^{n} f_k H^T(f_k)\|_i \to 0$ since $f_{k_i} \to 0$ as you get closer to $x^*$. This implies that $\xi_i \to 0$ as $x \to x^*$. Now, since $x^*$ has a stability region in the quotient gradient system, equation (3), the convergence of the QGS-based $\psi tc$ method to $x^*$ can be analyzed using the analysis in [33].

In summary, the QGS-based $\psi tc$ with $h_{max} = \infty$ will converge at least q-superlinerly if all the assumptions are statisfied, $\xi_i$ is sufficiently small, and $\xi_i \to 0$ as $x_i \to x^*$. The convergence of the QGS-based $\psi tc$ method to $x^*$ is locally q-quadratic if $\xi_i = O\|Q(x_i)\|$.

The solutions/equilibrium points of the QGS system (3) can be of three types:
1. A solution $x^*$ such that $F(x^*) = 0$, and $-DF(x^*)^T F(x^*) = 0$.
2. A solution $x^*$ such that $F(x^*) \neq 0$, and $-DF(x^*)^T F(x^*) = 0$.
3. A $x^*$ such that $F(x^*) = 0$, $-DF(x^*)^T F(x^*) = 0$, and $DF(x^*)^T$ is singular.

Unlike the standard $\psi tc$ or the Euler method, the QGS-based $\psi tc$ method when applied to (3) converges only to solutions of type 1 because of the QGS structure-dependent Jacobian approximation. This implies that approximating the QGS Jacobian guarantees that if the method converges, it will only converge to solutions of (1).

*Theorem 3: If the QGS-based $\psi$tc method converges to $x^*$, then $x^*$ is a solution of type 1 and hence, a solution of the original system (1).*

Proof:
If the QGS-based $\psi tc$ method converges, then it implies that $DF(x^*)^T F(x^*) = 0$ and the inverse $\left(h_{x^*}^{-1} I + DF(x^*)^T DF(x^*)\right)^{-1}$ exists where $h_{x^*}$ is the final value of the time step.
If $x^*$ is of type 3, then $DF(x^*)^T DF(x^*)$ is singular. Since the function $\lambda_{min}(DF(x)^T DF(x))$, the minimum eigenvalue of $DF(x)^T DF(x)$, is continuous in $x$, the sequence $\{\lambda_{i,min} = \lambda_{min}(DF(x_i)^T DF(x_i))\}$ converges to zero as $x_i \to x^*$. Thus, there exists a value of $i < \infty$ where $DF(x_i)^T DF(x_i)$ is singular and $h_i$ is very large such that $\left(h_i^{-1} I + DF(x_i)^T DF(x_i)\right) \approx DF(x_i)^T DF(x_i)$, which is a contradiction since $\left(DF(x_i)^T DF(x_i)\right)^{-1}$ must exist for convergence. Hence, the QGS-based $\psi tc$ method cannot converge to solutions of type 3.
If $x^*$ is of type 2, then the $Null(DF(x^*)^T) \neq \emptyset$, and $F(x^*) \in Null(DF(x^*)^T)$. This implies that $DF(x^*)^T DF(x^*)$ is singular, which implies that the QGS-based $\psi tc$ method cannot converge to solutions of type 2, since it will be a contradiction. Thus, if the QGS-based $\psi tc$ method converges, the solution can only be of type 1 and hence, the solution will be a solution of (1).

## V. NUMERICAL EXAMPLES

The numerical simulations were performed on a computer with an Intel® Core™ i7-3630QM CPU @2.40GHz processor and 16GB memory. All the simulations were performed with Matlab 7.11. Matlab's Fsolve is used as a check for accuracy. The systems simulated were the structure-preserving models of the WSCC 9-bus 3-machine system and the IEEE 145-bus 50-machine system with classical generators and constant impedance load models. The generalized list of equations for the structure-preserving model is as shown below.

Structure Preserving Model: For n generators and m buses,
$$\dot{\tilde{\delta}}_i = \tilde{\omega}_i \quad (19)$$
$$M_i \dot{\tilde{\omega}}_i = -D_i \tilde{\omega}_i + P_{m_i} - \frac{E'_{qi} V_i \sin(\tilde{\delta}_i - \tilde{\theta}_i)}{X'_{di}} - \frac{M_i}{M_i} P_{COI} \quad (20)$$

For generator buses $i = 1, \dots n$:
$$(I_{di} + jI_{qi})e^{-j(\delta_i - \pi/2)} = \sum_{k=1}^{m} Y_{ik} e^{j\alpha_{ik}} V_k e^{j\tilde{\theta}_k}$$
$$I_{di} = \frac{E'_{qi} - V_i \cos(\tilde{\delta}_i - \tilde{\theta}_i)}{X'_{di}}, \quad I_{qi} = \frac{V_i \sin(\tilde{\delta}_i - \tilde{\theta}_i)}{X'_{qi}} \quad (21)$$

For load buses $i = n + 1, \dots m$:
$$0 = \sum_{k=1}^{m} Y_{ik} e^{j\alpha_{ik}} V_k e^{j\tilde{\theta}_k} \quad (22)$$

$$\delta_0 = \frac{1}{M_T} \sum_{i=1}^{n} M_i \delta_i, \quad \omega_0 = \frac{1}{M_T} \sum_{i=1}^{n} M_i \omega_i$$
$$M_T = \sum_{i=1}^{n} M_i, \tilde{\delta}_i = \delta_i - \delta_0, \tilde{\omega}_i = \omega_i - \omega_0,$$
$$\tilde{\theta}_i = \theta_i - \theta_0 \text{ for } i = 1, \dots n,$$
$$P_{COI} = \sum_{i=1}^{n} P_{m_i} - \sum_{i=1}^{n} \frac{E'_{qi} V_i \sin(\tilde{\delta}_i - \tilde{\theta}_i)}{X'_{di}}$$

where $\delta_i$, $\omega_i$, $M_i$, $D_i$, $P_{m_i}$, $E'_{qi}$, $X'_{di}$, $V_i$, $\theta_i$ and $Y_{ik} e^{j\alpha_{ik}}$ are: rotor angle of machine $i$, speed of machine $i$, moment of inertia of machine $i$, damping of machine $i$, mechanical power of machine $i$, equivalent transient quadrature internal voltage of machine $i$, equivalent direct transient reactance of machine $i$, voltage magnitude at bus $i$, voltage angle at bus $i$, and the network admittance between buses $i$ and $k$, respectively.

TABLE I
CONTINGENCY LIST OF WSCC 9-BUS 3-MACHINE SYSTEM

| Contingency Number | Fault Bus | From Bus | To Bus |
|---|---|---|---|
| 1 | 7 | 7 | 5 |
| 2 | 7 | 8 | 7 |
| 3 | 6 | 4 | 6 |
| 4 | 6 | 6 | 9 |
| 5 | 9 | 9 | 8 |

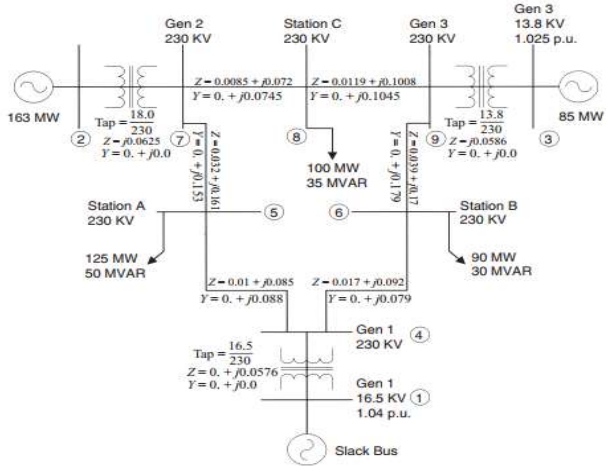

Fig. 2. The WSCC 9-bus 3-machine system. The value of Y is half the line charging.

The size and connectedness of the convergence region of the QGS-based $\psi tc$ method when used for UEP computation will be one of the metrics used in the evaluation of the proposed method's performance in comparison to other solvers. The convergence region of a UEP of a dynamic system for a numerical solver $N$ is defined as the set of initial points that converge to the UEP for the numerical solver [1]. Obviously, a large and connected convergence region implies that the initial point for the UEP computation does not have to be very close to the UEP for the algebraic solver to converge.

The convergence region of a UEP for a classical generator model can be constructed using a dimension-reduced model in which $\omega$ is a zero vector. This construction is made by creating a grid of initial points around the UEP in the machine angle space less the reference machine angle variable. The reference machine angle and the algebraic variables corresponding to the initial points on the grid are then updated by using the COI equation for machine angles and the solutions of the network equations at the grid points, respectively. The solutions of the equilibrium equations for the dynamic system starting at these initial points is then computed using the algebraic solver for which a convergence region of a UEP is being constructed. If the L2 norm of the difference between the computed equilibrium point and the UEP is below a defined threshold, then the initial point is in the convergence region of the UEP for that algebraic solver.

### A. The WSCC 9-Bus 3-Machine System with a Classical Generator Model

The method is tested on the WSCC 9-bus 3-machine system [1] (see Fig. 2), to compute a UEP on the stability boundary of a post-fault system. A uniform damping of $\lambda = 0.1$ is assumed, and the simulation is done in the Center of the Inertia (COI) reference framework. Initial time step $h_0$ of the $\psi tc$ and QGS-based $\psi tc$ methods are set to 0.1. Table I shows the list of contingencies used in our simulations of this system.

We first look at a simulation example where the initial point for a UEP computation is outside the convergence region of the NR method but within the convergence region of the proposed method. Fig. 3 shows a comparison of projected convergence regions of the NR method and the proposed method for the controlling UEP of contingency 1. We observe that an initial point, depicted by the black asterisk falls within the convergence region of the proposed method but outside the convergence region of the NR method. Hence, the proposed method can converge to the CUEP of contingency 1 while the NR method diverges. Fig. 3 also show that, the convergence region of the proposed method is the largest one among the two methods, and is much more compact and connected while that of the NR and the Continuous NR method is disconnected. The connected section of the convergence region of the QGS-based $\psi tc$ method is also larger. Consequently, the initial point for a UEP computation could be close and the NR method will still diverge, a less likely case with the proposed method.

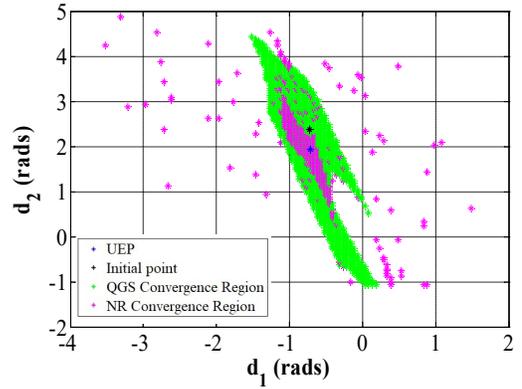

Fig. 3. Convergence regions of a UEP for the structure-preserving model (contingency 1) of a post-fault WSCC 9-bus 3-machine system. Convergence region of the NR method superimposed on the convergence region of the QGS-based $\psi tc$ method, with an example of an initial point that is in the convergence region of the QGS-based $\psi tc$ for a UEP but not in the convergence region of the NR method.

TABLE II
COMPARISON OF CONVERGENCE REGIONS OF METHODS FOR SOLVING THE CUEPs OF THE STRUCTURE -PRESERVING MODEL OF THE WSCC 9-BUS 3-MACHINE SYSTEM $\epsilon = 10^{-6}$

| Contingency | Initial Points in Connected Convergence Region | | | Initial Points Outside Connected Convergence Region | | |
|---|---|---|---|---|---|---|
| | NR | QGS $\psi tc$ | CNR | NR | QGS $\psi tc$ | CNR |
| 1 | 241 | 1016 | 426 | 62 | 0 | 160 |
| 2 | 231 | 337 | 293 | 26 | 1 | 24 |
| 3 | 230 | 490 | 346 | 120 | 3 | 193 |
| 4 | 215 | 510 | 345 | 127 | 1 | 221 |
| 5 | 240 | 361 | 324 | 68 | 0 | 26 |

Table II shows a comparison of the convergence regions of the controlling UEPs for the structure-preserving model of the post-fault WSCC 9-bus 3-machine system for the five contingencies in Table I using the NR method, the QGS-based $\psi tc$, and the continuous Newton-Raphson method (CNR). We observe that very few of the initial points that converge to a UEP using the QGS-based $\psi tc$ method, less than 1%, fall outside the connected portion of the convergence region of the UEP. On the other hand, 10% to 37% and 7% to 39% of the initial points that converge to a UEP using the NR and the CNR method fall outside the connected region of their convergence regions, respectively. This result show that larger portions of the UEP convergence region of the QGS-based $\psi tc$ method are connected compared the NR and the CNR methods. It also shows that the QGS-based $\psi tc$ method has a larger UEP convergence region compared to the NR method for the 5 contingencies. We also observe the UEP convergence region for the QGS-based $\psi tc$ method is larger than that of the CNR method for 3 contingencies. However, for the two contingencies where the UEP convergence region for the CNR method is larger than that for the QGS-based $\psi tc$ method, a third of the convergence region of the CNR method is not part of the connected portion, and consequently the connected portion of the convergence region of the QGS-based $\psi tc$ method is larger than the convergence region of the CNR method.

TABLE III
COMPARISON OF METHODS FOR SOLVING THE CUEPs OF THE STRUCTURE - PRESERVING MODEL OF THE WSCC 9-BUS 3-MACHINE SYSTEM $\epsilon = 10^{-6}$

| Method | Average Iterations | Average Computation Time (secs) |
|---|---|---|
| NR | 4 | 0.023 |
| Matlab's Fsolve | 5 | 0.035 |
| Original $\psi tc$ | 9 | 0.200 |
| QGS $\psi tc$ | 9 | 0.027 |
| QGS $\psi tc$ with NR | 8 | 0.027 |
| Continuous Newton | 15 | 0.051 |

Table III shows a comparison of the average number of iterations and computation time in seconds for CUEP evaluations for the 5 contingencies given an initial point, using Matlab's Fsolve, the NR method, the CNR [30] method, the original $\psi tc$ method (thus, the $\psi tc$ method applied to (3) or (10)), the QGS-based $\psi tc$, and the QGS-based $\psi tc$ combined with the NR method. The table shows that, by using the QGS-based $\psi tc$ method, we maintain the average number of iterations while improving the UEP computation speed by 7.4 times on average compared to the original $\psi tc$ method. We also get a lesser number of average iterations when the QGS-based $\psi tc$ method is combined with the NR method. The table also shows that the computation speed of the QGS-based $\psi tc$ method is comparable to computation speed of the NR method and Matlab's Fsolve and about 2 times faster compared to the CNR method for this numerical example. As expected, the QGS-based $\psi tc$ requires more iterations than both Matlab's Fsolve and the NR method. The CNR method requires the largest number of iterations for all the contingencies in this study.

TABLE IV
CONTINGENCY LIST OF IEEE 145-BUS 50-MACHINE SYSTEM

| Contingency Number | Fault Bus | From Bus | To Bus |
|---|---|---|---|
| 1 | 6 | 7 | 6 |
| 2 | 72 | 59 | 72 |
| 3 | 116 | 115 | 116 |
| 4 | 100 | 100 | 72 |
| 5 | 91 | 91 | 75 |
| 6 | 112 | 112 | 69 |
| 7 | 101 | 101 | 73 |
| 8 | 6 | 6 | 1 |
| 9 | 59 | 59 | 103 |

B. *The IEEE 145-Bus 50- Generator System with the Classical Generator Model*

The proposed method is also tested on the IEEE 145-bus 50-machine system for CUEP computations. A uniform damping of $\lambda = 0.5$ is assumed and the simulation is done in the Center of Inertia (COI) reference framework. An initial time step, the same as in the previous study, is used. Table IV shows the list of 9 contingencies used in the test.

TABLE V
COMPARISON OF METHODS FOR SOLVING THE CUEPs OF THE STRUCTURE - PRESERVING MODEL OF THE IEEE 145-BUS 50-MACHINE SYSTEM $\epsilon = 10^{-6}$

| Method | Average Iterations | Average Computation Time (secs) |
|---|---|---|
| NR | 3 | 0.038 |
| Matlab's Fsolve | 4 | 0.068 |
| Original $\psi tc$ | 9 | 10.411 |
| QGS $\psi tc$ | 9 | 0.078 |
| QGS $\psi tc$ with NR | 7 | 0.066 |
| Continuous Newton | 14 | 0.293 |

Table V shows a comparison of the average number of iterations and computation time in seconds for CUEP evaluations for the 9 contingencies given an initial point, using Matlab's Fsolve, the NR method, the CNR [30] method, the original $\psi tc$ method, the QGS-based $\psi tc$, and the QGS-based $\psi tc$ combined with the NR method. The table shows that, by using the QGS-based $\psi tc$ method, we maintain the number of iterations while improving the UEP computation speed by 133 times on average compared to the original $\psi tc$ method. We get a better number of average iterations, and a UEP computation speed of 157 times the speed of the original $\psi tc$ method when the QGS-based $\psi tc$ method is combined with the NR method. The table shows that the QGS-based $\psi tc$ method is about 4 times faster compared to the CNR method, and 2 times slower than the NR method. The table also shows that the computation speed of the QGS-based $\psi tc$ method is comparable to computation speed of Matlab's Fsolve for this numerical example. As expected, the QGS-based $\psi tc$ requires more iterations than both Matlab's Fsolve and the NR method. The CNR method requires the greatest number of iterations for all the contingencies in this study too.

C. *Discussion*

The results from the numerical simulations show that the proposed method introduces robustness into the computation of UEPs by expanding the connected portion of the convergence region, at the expense of a decrease in speed. Proposition 1 show that the transformation of (1) to (3) changes a UEP to an SEP and, in effect, creates a stability region making the application of TJM methods possible. From the results in Table III & V we observe that exploiting the structure of the QGS (3) by using an appropriate approximation of the Jacobian (11) can improve the speed of the proposed method significantly and ensure convergence only to solutions of or the original problem (1). The results in Fig. 3 and Table II show that the proposed method makes it possible for the initial guess for a UEP to be further away from the UEP compared to the NR method or the CNR method because the method has a larger connected convergence region.

As the QGS-based $\psi tc$ method is slower than the NR method, we recommend that the QGS-based $\psi tc$ method be used as a re-starting algorithm only when the NR method fails.

## VI. Conclusion

In this paper, we have proposed a new method that combines a QGS transformation with a TJU method for the computation of a UEP, given an initial point. The method converts the UEP to an SEP by changing the problem into a quotient gradient system. It then applies a quasi-Newton form of the pseudo-transient continuation method by exploiting the structure of the proposed QGS's Jacobian. The main advantages of this method are: 1) the method has a larger continuous convergence region than the NR method and thus, the initial guess does not have to be as close to the equilibrium point as the NR method requires; 2) it is faster than simply applying the exact pseudo-transient continuation method to the QGS; and 3) the proposed inexact pseudo-transient continuation method can only converge to solutions of the original system, unlike the exact pseudo-transient continuation method. The applicability of the proposed method to improving the robustness of CUEP computations in the direct method for transient stability analysis of power systems will be further investigated.

**Robert Owusu-Mireku** received his B.Sc degree in electrical and electronic engineering from the Kwame Nkrumah University of Science and Technology, Kumasi, Ghana and his M.Eng and M.S in electrical engineering from Cornell University, Ithaca, NY, USA. He is currently pursuing a Ph.D. degree at the School of Electrical and Computer Engineering in Cornell University, Ithaca, NY, USA. His current research interests include nonlinear system theory, nonlinear computation, renewable energy, and their practical applications to electric power systems.

**Hsiao-Dong Chiang** (M'87-SM'91-F'97) received his Ph.D degree in electrical engineering and computer science from the University of California Berkeley, Berkeley, CA, USA. He is a professor in the School of Electrical and Computer Engineering at Cornell University, Ithaca, NY, USA. He and his research team have published more than 350 referred papers. His current research interests include nonlinear system theory, nonlinear computation, nonlinear optimization, and their practical applications. He holds 17 U.S. and overseas patents and several consultant positions. He is the author of two books: *Direct Methods for Stability Analysis of Electric Power Systems: Theoretical Foundation, BCU Methodologies, and Application* (New York, NY, USA: Wiley, 2011) and (with Luis F. Alberto) Stability Regions of Nonlinear Dynamical Systems: Theory, Estimation, and Applications (Cambridge Univ. Press, 2015). He has served as an associate editor for serveral IEEE transactions and journals and is the founder of Bigwood Systems, Inc. in Ithaca, NY 14850, USA.